%% file: NewtonFBP2013R1.tex
\documentclass[12pt,a4paper]{article}
\usepackage[T1]{fontenc}
\usepackage{amsmath,amssymb,amsfonts}        % moduli per le estensioni ed i fonts AMS
\usepackage{mathptmx}
\usepackage[scaled=0.91]{helvet}
\usepackage{courier}
\usepackage[dvips]{graphicx}
\usepackage{psfrag}
\usepackage{color}
\usepackage[pdfauthor={Riccardo Fazio},pdftitle={Newton's Free Boundary Problem}, bookmarks,colorlinks]{hyperref}

\def\ds{\displaystyle}
\def\FBP{free boundary problem}

\def\D{\ensuremath{\mbox{D}}}

\title{A Non-Iterative Transformation Method for Newton's Free Boundary Problem}
\author{Riccardo Fazio \\
Department of Mathematics and Computer Science\\
University of Messina \\
Viale F. Stagno D'Alcontres, 31 \\
98166 Messina, Italy \\
E-mail: rfazio@unime.it \\
Home-page: \url{http://mat521.unime.it/fazio}}
\date{Submitted to International Journal Non-Linear Mechanics on September 3, 2013 and in revised form \today}
\begin{document}
\maketitle
\begin{abstract}
In book II of Newton's \textit{Principia Mathematica} of 1687 several applicative problems are introduced and solved.
There, we can find the formulation of the first calculus of variations problem that leads to the first \FBP \ of history.
The general calculus of variations problem is concerned with the optimal shape design for the motion of projectiles subject to air resistance.
Here, for Newton's optimal nose cone \FBP , we define a non-iterative initial value method which is referred in literature as transformation method.
To define this method we apply invariance properties of Newton's \FBP \ under a scaling group of point transformations.
Finally, we compare our non-iterative numerical results with those available in literature and obtained via an iterative shooting method.
We emphasize that our non-iterative method is faster than shooting or collocation methods and does not need any preliminary computation to test the target function as the iterative method or even provide any initial iterate.
Moreover, applying Buckingham Pi-Theorem we get the functional relation between the unknown free boundary and the nose cone radius and height.  
\end{abstract}
\smallskip

\smallskip


%\newpage

\section{Introduction}
%\textcolor{red}{Bigliograpy}
In book II of Newton's \textit{Principia Mathematica} of 1687 several applicative problems are introduced and solved.
There, we can find the first calculus of variations problem, predating the brachistochrone problem of the 1690s, that leads to the first \FBP \ of history.
The general calculus of variations problem is concerned with the optimal shape design for the motion of projectiles subject to air resistance, see Edwards \cite{Edwards:NNC:1997}.
\begin{figure}[!hbp]
\centering
\psfrag{x}[][]{$x$} 
\psfrag{y}[][]{$y$} 
\psfrag{P}[][]{$(r, h)$} 
\psfrag{r}[][]{$r$} 
\psfrag{h}[][]{$h$} 
\includegraphics[width=.6\textwidth,height=.6\textwidth]{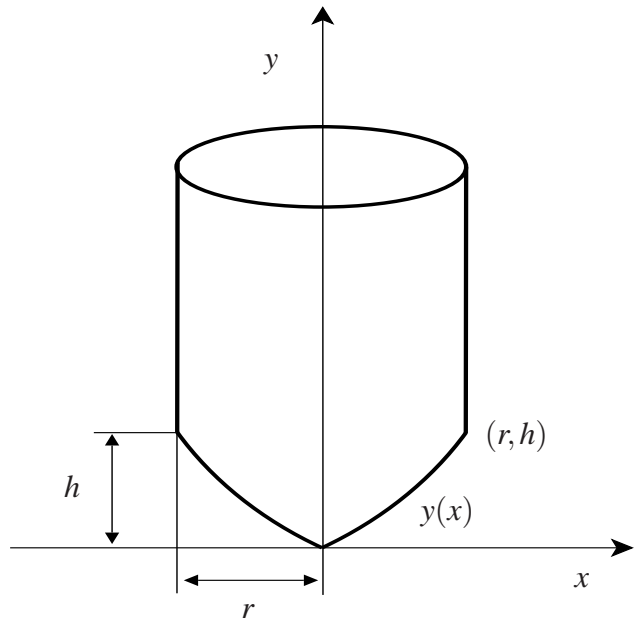}
\put(-80,40){$y(x)$}
\caption{Physical setup for Newton's projectile shape design.}
\label{fig:NewtonPro}
\end{figure}
Figure~\ref{fig:NewtonPro} shows a projectile shaped like a bullet with a nose cone having radius $r$ and height $h$; here $x$ and $y$ are Cartesian coordinates.
This nose cone is a surface of revolution determined by the plane curve $y = y(x)$. 
The term \lq \lq nose cone\rq \rq \ refers to the generic shape and not necessarily an actual \lq \lq cone\rq \rq .
On the basis of  reasonable assumptions on the air resistance, Newton established that the air resistance force $F$ acting to the projectile moving at a velocity $v$ is given by
\begin{equation}\label{eq:F}
F = 2 k\rho v^2 \ ,
\end{equation}
where $\rho$ is the density of the air, %$A = \pi r^2$ is the cross-sectional area of the cylindrical projectile, 
and the fundamental parameter $k$, indicated as the drag coefficient, is given by the formula
\begin{equation}\label{eq:k}
k = \int_0^r \frac{2 \pi x}{\left[\frac{dy}{dx}(x)\right]^2 + 1} \; dx \ .
\end{equation}
This model is based on the assumption of gas flows as independent movement of non-interacting mass particles, which hit the projectile shape and change their momentum. 
The absence of particle interactions is in contradiction to laminar and turbulent flow models, which are preferred for dense fluids. 
Thus, Newton model is only useful for three occasions: motion in low pressure gas, generating of good tasks for the calculus of variation and providing the first example of a \FBP .

Usually, given a specific configuration of the bullet researchers compute the reduced drag coefficient $k^* = k/\pi r^2$ instead of $k$, and in table~\ref{tab:k} we report the results concerning several nose cone shapes investigated by Newton.
\input{TAB1}
It is interesting to realize that a rounded hemisphere and a pointed cone provide the same value $k^* = 1/2$.
The paraboloid, with $y(1) = 1$, yields a smaller value.
For the conical frustum the function $y(x)$ is a straight segment forming an angle $\theta$ with the $x$-axis.
The optimally angled flat-tipped conical frustum offers least air resistance when $\tan 2 \theta = 2$.
Newton's optimal nose cone is defined by the solution of a \FBP , to be discussed and solved numerically in the next sections.
Newton's flat-tipped frustum offers less air resistance than all of the simple shapes.
There are two competing effects: the flat tip has large air resistance but it allows the nose cone to have steeper sides, which reduces the air resistance.
For over 300 years, Newton's solution stood as the minimizer but it is only the radially-symmetric flat-tipped minimizer.
A radially-symmetric nose cone with indented tip results in a smaller value of $k^*$, that is $k^* = 0.29519$, that Newton's optimal value, see Gallant \cite{Gallant:2012:DPS}. 

Landau \cite{Landau:HCM:1950} was the first to point out that \FBP s are always non-linear.
Therefore, this kind of problems are often solved numerically.
Moreover, normally \FBP s are transformed into boundary value problems (BVPs), see Ascher et al. \cite{Ascher:1981:RBV} or \cite[p. 471]{AscherBook}.
In this context, sometimes, it is possible to solve a given \FBP \ non-iteratively, see the survey by \cite{Fazio:1998:SAN}, whereas BVPs are usually solved iteratively.

Here, for Newton's \FBP \ optimal nose cone, we define a non-iterative initial value method which is referred in literature as non-iterative transformation method (ITM).
Indeed, non-ITMs can be defined within Lie's group invariance theory.
As far as group invariance theory is concerned, we refer to Bluman and Cole \cite{Bluman:1974:SMD}, Bluman and Kumei \cite{Bluman:1989:SDE}, Barenblatt \cite{Barenblatt:1996:SSI}, or Dresner \cite{Dresner:1999:ALT}.
The first application of a non-ITM to a \FBP \ was given by Fazio and Evans \cite{Fazio:1990:SNA}.
In the past, the main drawback of non-ITMs was that they were considered not widely applicable: see the critical considerations by Fox et al. \cite{Fox:1969:LBL}, Meyer \cite[pp. 97-98]{Meyer:1973:IVM}, Na\cite[p. 137]{Na:1979:CME} or Sachdev \cite[p. 218]{Sachdev:1991:NOD}. 
In fact, the simplest way in order to verify if a non-ITM is applicable to a particular problem is to use an inspectional analysis as shown by Seshadri and Na \cite[pp. 157-168]{Seshadri:1985:GIE} (cf. also the discussion on inspectional analysis by Birkhoff \cite[pp. 99-103]{Birkhoff:1960:HSL}).
In relation to the transformation of \FBP s to initial value problems (IVPs), it is also possible to define an iterative extension
of the TM which is always applicable \cite{Fazio:1994:FSEb,Fazio:1996:NAN,Fazio:2013:BPF}.

\section{Newton's \FBP}
\begin{figure}[!hbp]
\centering
\psfrag{x}[][]{$x$} 
\psfrag{x0}[][]{$a$} 
\psfrag{y}[][]{$y$} 
\psfrag{P}[][]{$(r, h)$} 
\includegraphics[width=.6\textwidth,height=.6\textwidth]{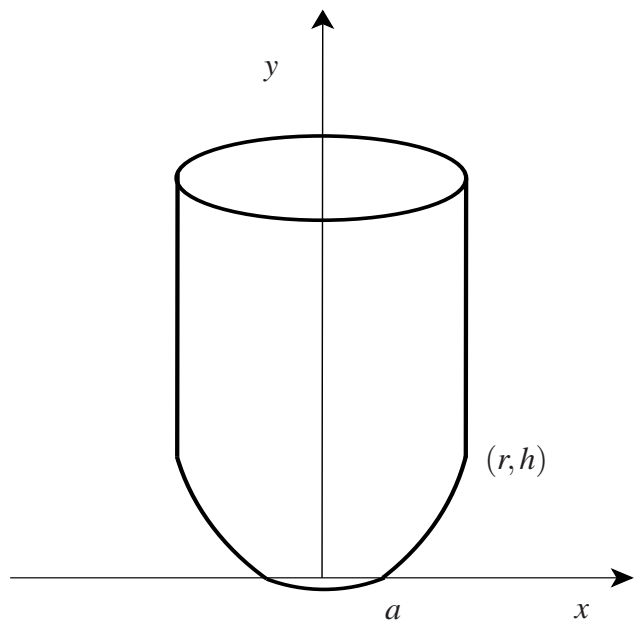}
\caption{Physical setup for Newton's optimal flat-tip projectile shape.}
\label{fig:NewtonPFBP}
\end{figure}

In \cite{Edwards:NNC:1997} Edwards exploits modern computer algebra tools to explore the origin and meaning of Newton's nose cone problem. 
In figure~\ref{fig:NewtonPFBP} we show the physical setup for Newton's optimal flat-tip projectile shape. 

Newton's intuition was that the optimal nose cone of least resistance would have a flat circular tip joined to the cylindrical body by a curvilinear band. 
In this particular case, the value of the reduced drag coefficient can be computed by
\begin{equation}\label{eq:k*}
k = a^2 + \int_a^r \frac{2 x}{\left[\frac{dy}{dx}(x)\right]^2 + 1} \; dx \ .
\end{equation}
But, what should be the radius $a$ of the flat tip and what should be the shape $y = y(x)$ of the arc generating this optimal band by revolution around the $y$-axis? 
Nowadays, we would regard this as a calculus of variations problem, complicated by a variable midpoint condition, and proceed to set up the appropriate Euler-Lagrange equation  
\begin{equation}
\left[\frac{\partial}{\partial y} - \frac{d}{dx}\frac{\partial}{\partial \dot{y}} \right] \Phi(x, y, \dot{y}) = 0 \ ,
\end{equation} 
where we have used Newton's dot notation for the first derivative and
\begin{equation}
\Phi(x, y, \dot{y}) = \frac{2 x}{\left(\dot{y}\right)^2 + 1} \ .
\end{equation}
Therefore, we get the free BVP
\begin{align}\label{Newton:FBP}
&{\ds \frac{d^2y}{dx^2}}
= \frac{\ds \frac{dy}{dx}\left[\left(\frac{dy}{dx}\right)^2+1\right]}{x\left[3\left(\ds \frac{dy}{dx}\right)^2-1\right]}\ , \qquad x \in [a, r] \nonumber \\[-.5ex]
& \\[-.5ex]
& y(a) = 0 \ , \quad {\ds \frac{dy}{dx}(a) = 1} \ , \qquad y(r) = h \ ,  \nonumber
\end{align}
%$r$ and $h$ are the radius and the nose height  of the projectile, respectively, 
where the boundary condition on the derivative at $x = a$ means that the tangent to the arc $y(x)$ at $a$ has the same direction of $y = x$, the other two boundary conditions comes from the geometrical configuration of the projectile, and $a$, the length of the frustum, is the unknown free boundary.
The second order ordinary differential problem (\ref{Newton:FBP}) is Newton's optimal \FBP .

\section{Scaling invariance and a non-ITM}
Let us define a non-ITM for the numerical solution of Newton's \FBP .
The governing differential equation and the two free boundary values in (\ref{Newton:FBP}) are invariant under the scaling group of point transformation
\begin{equation}\label{eq:sgroup}  
x^* = \lambda  x \ , \quad a^* = \lambda  a \ , \quad y^* = \lambda y  \ ,
\end{equation}
where $\lambda$ is the group parameter.
Moreover, the end-point condition in (\ref{Newton:FBP}), the condition at $x= r$, is not invariant with respect to (\ref{eq:sgroup}).
Using the invariance properties, we can define the following non-iterative algorithm for the numerical solution of Newton's optimal \FBP \ (\ref{Newton:FBP})
\begin{enumerate}
\item Input $a^* > 0$, $r > 0$ and $h > 0$. %and $\Delta x$.
%\item Define a uniform grid on $[0, 1]$ with step size $\Delta x$.  
\item Solve the auxiliary IVP
\begin{align}\label{eq:ivp}
&{\ds \frac{d^2y^*}{dx^{*2}}
= \frac{\ds \frac{dy^*}{dx^*}\left[\left(\frac{dy^*}{dx^*}\right)^2+1\right]}{x^*\left[3\left(\ds \frac{dy^*}{dx^*}\right)^2-1\right]}} \ , \qquad x \in [a^*, r^*]  \nonumber \\[-1.5ex]
& \\[-1.5ex]
& y^*(a^*) = 0 \ , \qquad {\ds \frac{dy^*}{dx^*}(a^*)} =  1 \ , \nonumber 
\end{align}
where $r^*$ is defined implicitly by the condition
\begin{equation}\label{eq:event}
r^* = r \, \frac{y^*(r^*)}{h} \ .
\end{equation}
\item Compute $\lambda$ by 
\begin{equation}\label{eq:lambda1}
\lambda = \frac{y^*(r^*)}{h} \ .
\end{equation}
\item Rescale the numerical solution to get $y(x)$ and $\frac{dy}{dx}(x)$ according to (\ref{eq:sgroup}),
that is
\begin{equation}
y(x) = \lambda^{-1} y^*(x^*) \ , \quad 
{\ds \frac{dy}{dx}(x)} =  {\ds \frac{dy^*}{dx^*}(x^*)} \ ,
\end{equation}
in particular we find that
\begin{equation}
{\ds \frac{dy}{dx}(r)} =  {\ds \frac{dy^*}{dx^*}(r^*)} \ ,  
\end{equation}
and rescale the free boundary
\begin{equation}\label{eq:FB2}
a = \lambda^{-1} a^* \ .
\end{equation}
\end{enumerate}
In order to apply the condition (\ref{eq:event}) we can use an {\it event locator}.
A simple event locator is defined below.
We integrate the auxiliary IVP (\ref{eq:ivp}) until we get at a mesh point $x_k^*$ where $x_k^* < r \, y^*_k/h$, %and we compute  
%\begin{equation}\label{eq:evloc1}
%x_{u}^* = x_k^* + (h \, r^*/r - y_{k-1}^*) \frac{\Delta x^*}{y_k^*-y_{k-1}^*} \ .
%\end{equation}
then, we repeat the last step with the smaller step size given by
\begin{equation}\label{eq:evloc2}
%\Delta x_u^* = x_{u}^* - x_k^* \ .
\Delta x_u^* = \left(h \, \frac{r^*}{r} - y_{k-1}^*\right) \frac{\Delta x^*}{y_k^*-y_{k-1}^*} \ .
\end{equation}
In defining the last step size in equation 
%(\ref{eq:evloc1})-
(\ref{eq:evloc2}), we apply a first order (linear) Taylor formula at $x_k^*$ where we have replaced the first derivative with a backward finite difference approximation. 
Let us remark here, that the more recent versions of MATLAB allow to define a problem related event locator since they have a build in event locator algorithm.

\subsection{Numerical results}
Let us consider first a specific problem and fix the values $r=h=1$ used by Newton and by Edwards \cite{Edwards:NNC:1997}.
For the numerical solution of the auxiliary IVP (\ref{eq:ivp}) we used a uniform grid, with $\Delta x^* = 0.001$ and $a^* =0.5$, and applied a second order Heun method (RK2) \cite{Heun:1900:BNI}, the classical fourth order Runge-Kutta method (RK4) \cite{Runge:1895:UNA,Kutta:1901:BNI}, and the sixth order Runge-Kutta method (RK6) reported by Butcher in \cite[p. 178]{Butcher:NMO:2003}.
All computations were performed in double precision arithmetic, for all methods the last step size $\Delta x_u^*$ was about $3.51\D{-04}$.
In table~\ref{tab:RKes1} we list sample numerical results. 
\input{TAB2}
Indeed, the numerical results obtained by RK4 and RK6 are very close.
More accurate numerical results might be obtained by a grid refinement.
%In table~\ref{tab:RKes1:Dx} we list the numerical results obtained by RK6 and a grid refinement with the reported step sizes. 
 %\input{TABes1Dx}

The obtained values can be compared with the values $a \approx 0.350942572$ and $\frac{dy}{dx}(1) \approx 1.916801246$ computed by the symbolic program Mathematica and reported by Edwards \cite{Edwards:NNC:1997}.
He applied a shooting method using {\bf NDSolve}, the numerical IVP solver of the Mathematica PSE, and the root-finder implemented in {\bf FindRoot} with initial iterates $a_0 = 0.5$ and $a_1 = 0.3$.
%In this way Edwards was able to compute the value $a \approx 0.350943$.
See also the discussion on the numerical solution of this problem by Shampine et al. \cite[pp. 205-206]{Shampine:SOM:2003}: they advice to rewrite this \FBP \ in standard form \cite{Ascher:1981:RBV} and to apply the MATLAB BVP solver {\bf bvp4c} with the initial iterate $y_0(x) = x$ and $a_0 = 0.5$.
For the sake of completeness, we have implemented a shooting method within MATLAB and followed Shampine and his co-authors suggestion.    
The results reported below were obtained by the IVP {\bf ODE45} solver and the BVP solver {\bf bvp4c},
from the MATLAB ODE suite written by Samphine and Reichelt \cite{Shampine:1997:MOS}, with the accuracy and adaptivity parameters defined by default.
As far as the shooting method is concerned, we implemented both the bisection and the secant root-finder, along with the simple termination criterion
\begin{equation}\label{eq:TC}
|2 (a_k-a_{k-1})| \le 1 \D{-06} |a_k+a_{k-1}| \ ,
\end{equation}
to find a zero of the implicit function
\begin{equation}\label{eq:shoot:F}
F(y(1; a)) = y(1; a) - 1 \ ,
\end{equation}
where $a$ has been indicated as a parameter for the IVP solver.
The condition (\ref{eq:TC}) was meet by the bisection method after 22 iterations and by the secant method after 7 iterations; both root-finders got as results $a \approx 0.350943$ and $\frac{dy}{dx}(1) \approx 1.916801$.
As far as the collocation method implemented in {\bf bvp4c} is concerned, we used 101 and 201 grid points and the software took 2121 and 4221 calls to the ordinary differential equation function plus 36 and, again, 36 calls to the boundary condition function, respectively, and in both cases got as results $a \approx 0.3509426$ and $\frac{dy}{dx}(1) \approx 1.9167981$. 
The numerical results, reported in table \ref{tab:RKes1}, show that we are able to get an accurate numerical approximation of the solution of (\ref{Newton:FBP}) non-iteratively.
We remark that our non-ITM is at least seven times faster the the shooting method used by Edwards \cite{Edwards:NNC:1997}.
Moreover, the selected starting values bracket the zero of the implicit function (\ref{eq:shoot:F}), one of them, namely $a_0 = 0.5$ is an obvious choice, but why not take also $a_1 = 0.25$ or $a_1 = 0.75$?
As it is easily seen, when using a shooting method, a few trial computations should be considered in order to get suitable values of the starting values and usually are worthwhile.
On the other hand, by using the non-ITM we can forget about choosing any starting iterate.
Let us notice that, for the application of the non-ITM, we were not able to consider the value $a^*=0$ because the governing differential equation of the auxiliary IVP (\ref{eq:ivp}) has a singularity at $x^* = 0$.
In any way, the proposed method is independent on the choice of $a^*$, and, in fact, comparable numerical results have been obtained by setting $a^* = 1$.

\begin{figure}[!hbp]
\centering
\psfrag{x}[][]{$x$} 
\psfrag{u}[][]{} 
%\psfrag{u}[][]{$y(x)$, $\frac{dy}{dx}(x)$} 
\psfrag{x*}[][]{$x^*$} 
\psfrag{u*}[][]{} 
%\psfrag{u*}[][]{$y^*(x^*)$, $\frac{dy^*}{dx^*}(x^*)$} 
\psfrag{y}[][]{$y$} 
\psfrag{dy}[][]{$\frac{dy}{dx}$} 
\psfrag{y*}[][]{$y^*$} 
\psfrag{dy*}[][]{$\frac{dy^*}{dx^*}$} 
\includegraphics[width=180pt,height=360pt]{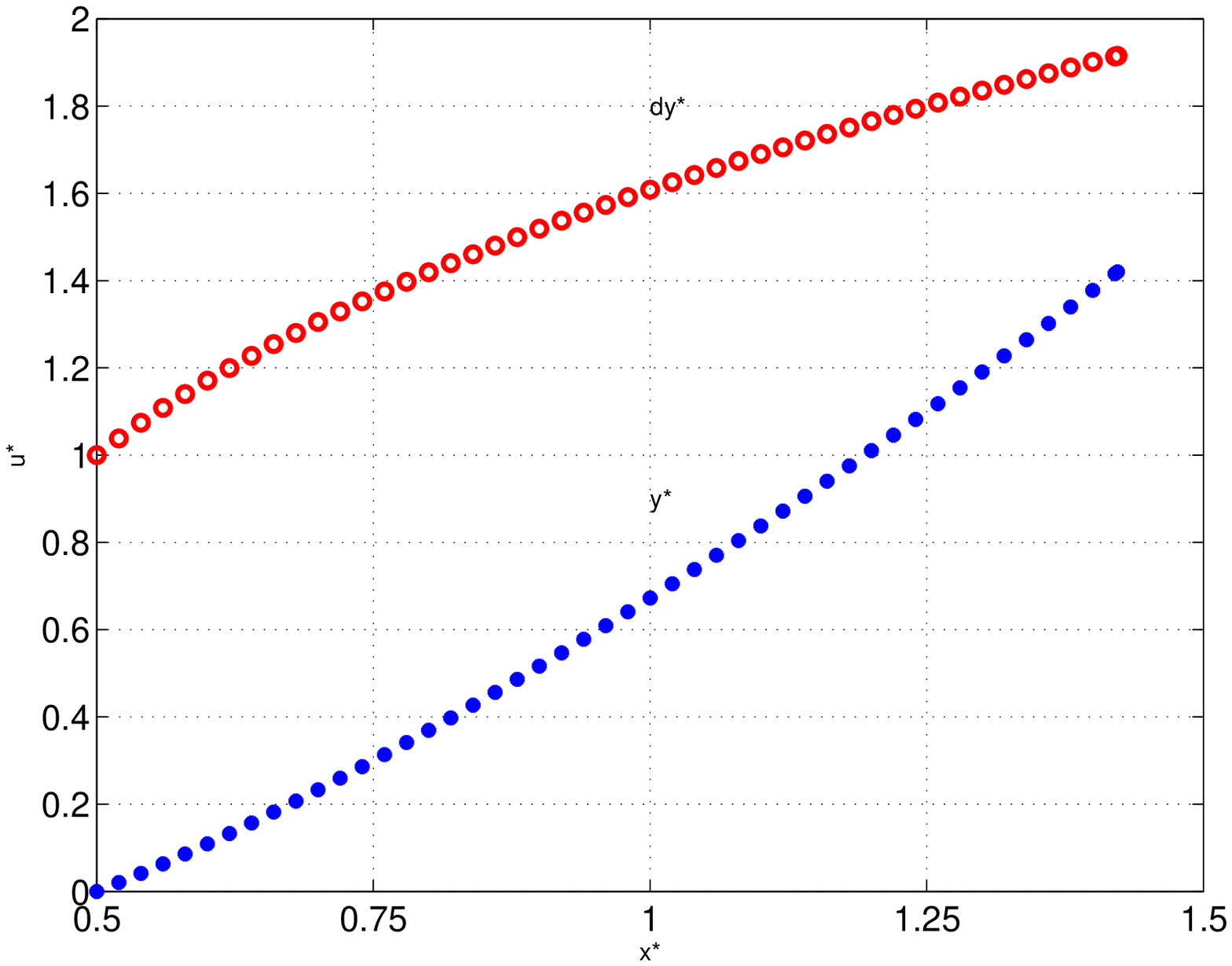} \hfil
\includegraphics[width=180pt,height=360pt]{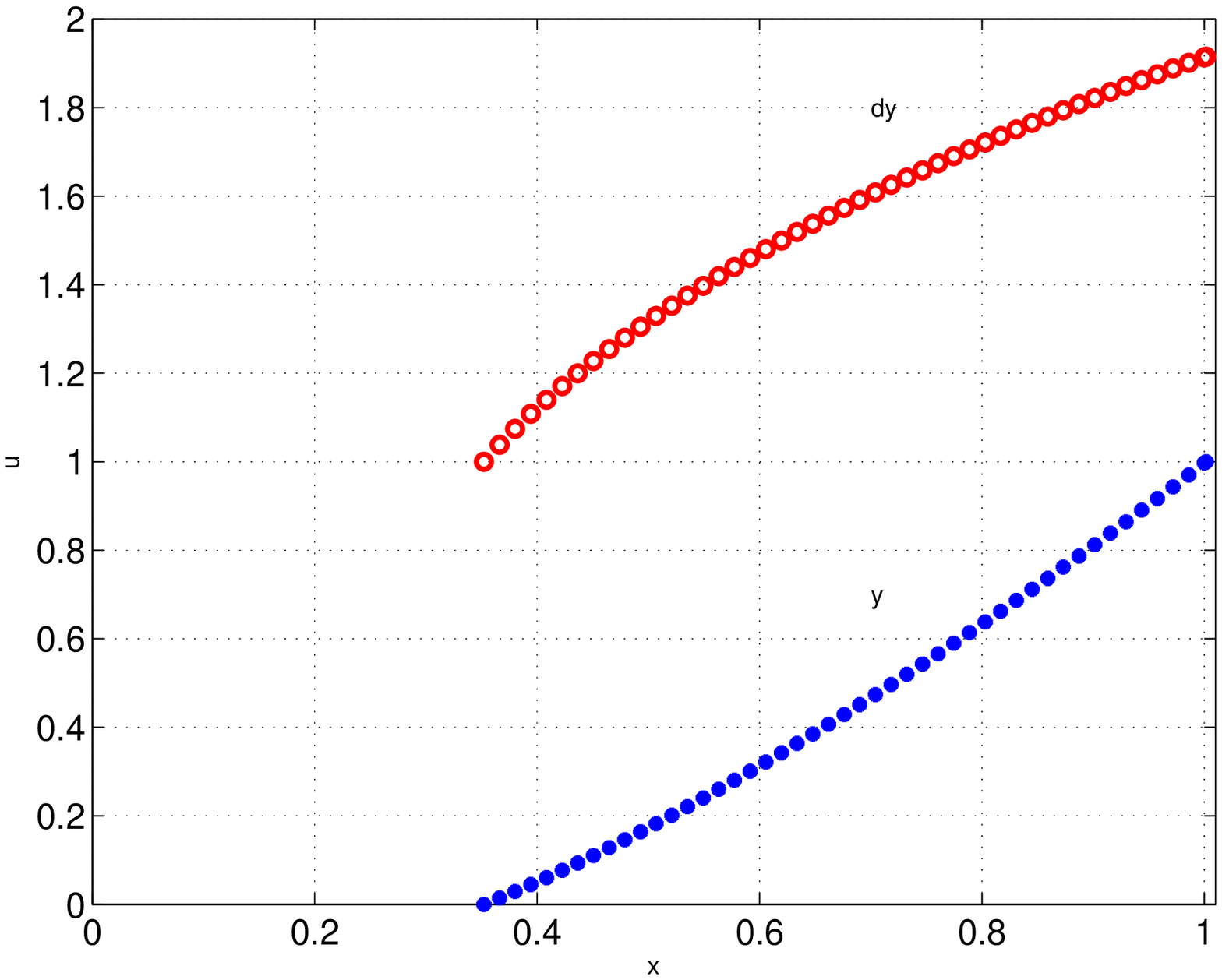}
\caption{On the left frame: solution of the initial value problem (\ref{eq:ivp}). On the right frame: solution of Newton's \FBP \ (\ref{Newton:FBP}). }
%The symbols stand for: $\triangledown$ $\frac{du_\e^*}{dx^*}(x^*)$, $\vartriangle$ $u_\e^*(x^*)$, $\square$ $\frac{du_\e}{dx}(x)$, and $\circ$ $u_\e(x)$.}
\label{fig:Es}
\end{figure}
In figure~\ref{fig:Es} we plot the numerical solution of the initial value problem (\ref{eq:ivp}), obtained by RK6 with $\Delta x^* = 0.02$, as well as the rescaled solution for the \FBP \ (\ref{Newton:FBP}).
From figure \ref{fig:Es} it is hardly possible to notice how the last step, $\Delta x_u^* \approx 0.0023$, is smaller with respect to the previous ones.

As far as different values of $r$ and $h$ are concerned, it is evident that $a = a(r,h)$ and applying Buckingham Pi-Theorem \cite{Buckingham:1914:PSS,Buckingham:1915:MEF} we get the functional form 
\begin{equation}\label{eq:f}
a = h f(r/h) \ .
\end{equation}
In fact, we have here three length quantities but only one independent fundamental dimensional unit so that we can introduce two a-dimensional quantities: $a/h$ and $r/h$.
Moreover, we notice that (\ref{eq:f}) has to be invariant with respect to the considered scaling invariance.
In fact, under the scaling group (\ref{eq:sgroup}) $r$ and $h$ transform as follows
\begin{equation}
r* = \lambda r \ , \qquad h* = \lambda h \ ,
\end{equation}
and, therefore,
\begin{equation}
f\left(\frac{r^*}{h^*}\right) = \frac{a^*}{h^*} =  \frac{a}{h} = f\left(\frac{r}{h}\right) \ .
\end{equation}
This means that: if the value of $f$ is known for a fixed value of the ratio $r/h$, then it is a simple matter to compute the corresponding value of $a$ by multiplying $h$ by this value of $f$. 

\begin{figure}[!htb]
\centering
\psfrag{f}[][]{$f(r)$} 
%\psfrag{x0}[][]{$a$} 
\psfrag{r}[][]{$r$} 
\includegraphics[width=.8\textwidth,height=7cm]{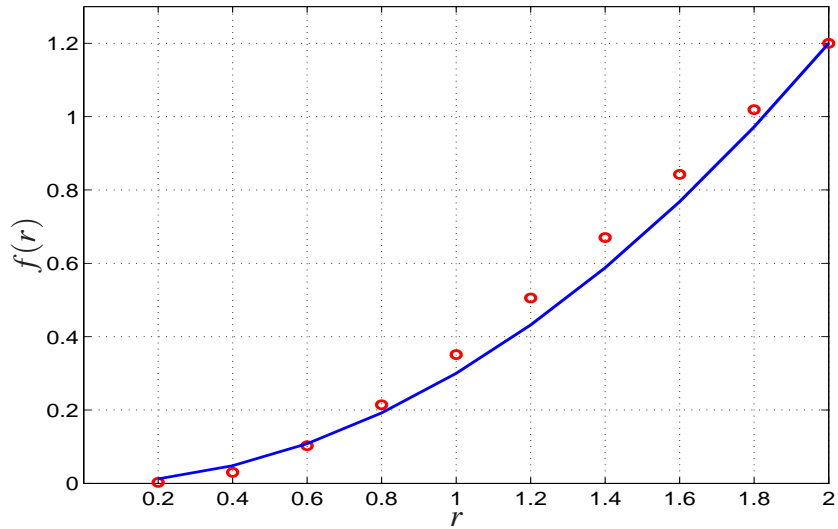}
\caption{The function $f(r/h)$ for $h=1$.}
\label{fig:fdir}
\end{figure}
Without lost of generality, we can take a length special unit defined by setting $h = 1$.
In this way we get the simpler relation $a = f(r)$.
Figure~\ref{fig:fdir} shows the behaviour of this function along with the parabola $ 0.3 r^2$ plotted for comparison. 
The coefficient $0.3$ is approximately the value that allows the last computed point, corresponding to $r = 2$, to belong to the parabola. 
It is surprising to recognize that also the first point, for $r = 0.2$ seems to belong approximately to the same parabola.
The results plotted in figure~\ref{fig:fdir} are in agreement with our intuition because as $r$ goes to zero we would expect that also $a$ goes to zero.

\section{Concluding remarks}
On account of the scaling invariance (\ref{eq:sgroup}), and by using the two invariants
\begin{equation}\label{eq:invariants}
u = \frac{y}{x} \ , \qquad s = \frac{dy}{dx} \ ,
\end{equation}   
we can reduce Newton's differential equation to a first order one, namely
\begin{equation}\label{eq:1or}
\frac{du}{ds} = \frac{(s-u)(3 s^2-1)}{s (s^2+1)} \ .
\end{equation}   
It turns out that equation (\ref{eq:1or}) can be solved in closed form
\begin{equation}\label{eq:1or:exact}
u = \frac{s (s^2 - \log(v) + 3 s^4/4) + C s}{(s^2+1)^2} \ ,
\end{equation}
where $C$ is an arbitrary constant.
If we substitute the invariants expressions (\ref{eq:invariants}) into (\ref{eq:1or:exact}), then we get a further first order differential equation, in the $x$ and $y(x)$ variables, whose solution is the general solution of the original second order differential equation.
Unfortunately, to solve this last first order differential equation seems to be a formidable task.
However, as we have seen in the previous section, the scaling invariance of a mathematical model can be used also to provide a simple way to get a numerical solution.

In the above context, in this paper we have defined a non-ITM for the numerical solution of Newton's \FBP \ (\ref{Newton:FBP}).
The obtained numerical result clearly indicate that this is an accurate and efficient way to deal with a problem of this kind. 
Naturally, the proposed non-ITM can be applied also to other \FBP s as discussed in \cite{Fazio:1990:SNA,Fazio:1998:SAN}.
In particular, the class of \FBP s 
\begin{align}\label{eq:FBPclass}
& {\displaystyle \frac{d^{2}y}{dx^2}} = y^{1-2\delta} \, 
\phi \left( x\, y^{-\delta},
{\displaystyle \frac{dy}{dx}}\, y^{\delta-1}
\right) \nonumber \\[-1ex]
& \\[-1ex]
& y(a) = \alpha \, a^{1/\delta} \ , \quad {\displaystyle \frac{dy}
{dx}} (a)  = \beta \, a^{(1/\delta)-1} \ ,  
\quad y(r) = h \ , \nonumber
\end{align}
where $\alpha$ and $\beta$ are given constants, can be characterized by the scaling group of transformations
\begin{equation}\label{eq:gscaleinv}
x^* = \lambda^{\delta} x \ , \quad a^* = \lambda^{\delta} a \ , \quad y^* = \lambda y  \ ,   
\end{equation}
where $\delta$ is a non-zero parameter.
It is a simple matter to verify that Newton's \FBP \ (\ref{Newton:FBP}) belongs to (\ref{eq:FBPclass}).
Of course, different classes of \FBP s can be characterized by imposing the invariance with respect to other groups of point transformations, e.g. the spiral \cite{Fazio:2010:TGN} or the translation group \cite{Fazio:2010:ESI}.

\bigskip
\bigskip

\noindent
{\bf Acknowledgements.} 
%The author is grateful to an anonymous reviewer .
This work was partially supported by GNCS of INDAM and by the University of Messina.
%\clearpage
%\newpage

\end{document}

%% file: TAB1.tex
\begin{table}[!hbt]
\caption{Different nose cone shapes and the corresponding values of the reduced drug coefficient ${k^*}$ computed for $r = h = 1$.}%, by equation \ref{eq:k}.}
\vspace{.5cm}
\renewcommand\arraystretch{1.3}
	\centering
		\begin{tabular}{lr@{.}l}%r@{.}l}%r@{.}l}
\hline 
Shape
%& \multicolumn{2}{c}%
%{$\e = 1.0\mbox{D}-03$}
%& \multicolumn{2}{c}%
%{${\ds \frac{du_\e}{dx}(0)}$}
& \multicolumn{2}{c}%
{${k^*}$} \\[1.2ex]
\hline
Hemisphere &  0 & 5000 \\  
Pointed cone &  0 & 5000 \\  
Paraboloid &  0 & 4024 \\  
Optimal conical frustum &  0 & 3820 \\  
Newton's optimal nose cone &  0 & 3748 \\  
\hline			
		\end{tabular}
	\label{tab:k}
\end{table}

%% file: TAB2.tex
\begin{table}[!hbt]
\caption{Sample numerical results for (\ref{Newton:FBP}) with $r = h = 1$.}%1.0\mbox{D}-03$. Here and in the following the $\mbox{D}-k = 10^{-k}$ means a double precision arithmetic.}
\vspace{.5cm}
\renewcommand\arraystretch{1.3}
	\centering
		\begin{tabular}{cr@{.}lr@{.}l}%r@{.}l}
\hline 
Method 
%& \multicolumn{2}{c}%
%{$\e = 1.0\mbox{D}-03$}
& \multicolumn{2}{c}%
{${\ds \frac{dy}{dx}(r)}$}
& \multicolumn{2}{c}%
{$x_0$} \\[1.2ex]
\hline
RK2 &  %13 & 8208 &     
1 & 916561011522854 & 0 & 351123613134112 \\
RK4 &  %13 & 8152 &     
1 & 916560741682499 & 0 & 351123613136137 \\
RK6 &  %13 & 8152 &     
1 & 916560741682204 & 0 & 351123613136137 \\
\hline			
		\end{tabular}
	\label{tab:RKes1}
\end{table}

%% file: NewtonFBP2013R1.bbl
\begin{thebibliography}{10}

\bibitem{AscherBook}
U.~M. Ascher, R.~M.~M. Mattheij, and R.~D. Russell.
\newblock {\em Numerical Solution of Boundary Value Problems for Ordinary
  Differential Equations}.
\newblock Prentice Hall, Englewood Cliffs, New Jersey, 1988.

\bibitem{Ascher:1981:RBV}
U.~M. Ascher and R.~D. Russell.
\newblock Reformulation of boundary value problems into \lq \lq standard\rq \rq
  \ form.
\newblock {\em SIAM Rev.}, 23:238--254, 1981.

\bibitem{Barenblatt:1996:SSI}
G.~I. Barenblatt.
\newblock {\em Scaling, Self-Similarity and Intermediate Asymptotics}.
\newblock Cambridge University Press, Cambridge, 1996.

\bibitem{Birkhoff:1960:HSL}
G.~Birkhoff.
\newblock {\em Hydrodynamics: A Study on Logic, Fact and Similitude}.
\newblock Princeton University Press, 1950, Princeton, 2nd edition, 1960.

\bibitem{Bluman:1974:SMD}
G.~W. Bluman and J.~D. Cole.
\newblock {\em Similarity Methods for Differential Equations}.
\newblock Springer, Berlin, 1974.

\bibitem{Bluman:1989:SDE}
G.~W. Bluman and S.~Kumei.
\newblock {\em Symmetries and Differential Equations}.
\newblock Springer, Berlin, 1989.

\bibitem{Buckingham:1914:PSS}
E.~Buckingham.
\newblock On physically similar systems.
\newblock {\em Physical Rev.}, 4:354--376, 1914.

\bibitem{Buckingham:1915:MEF}
E.~Buckingham.
\newblock Model experiments and the forms of empirical equations.
\newblock {\em Trans. Am. Mech. Eng.}, 37:263--296, 1915.

\bibitem{Butcher:NMO:2003}
J.~C. Butcher.
\newblock {\em Numerical Methods for Ordinary Differential Equations}.
\newblock Whiley, Chichester, 2003.

\bibitem{Dresner:1999:ALT}
L.~Dresner.
\newblock {\em Applications of {L}ie's Theory of Ordinary and Partial
  Differential Equations}.
\newblock Institute of Physics Publishing, London, 1999.

\bibitem{Edwards:NNC:1997}
C.H. Edwards.
\newblock Newton's nose-cone problem.
\newblock {\em Mathematica J.}, 7:64--71, 1997.

\bibitem{Fazio:1994:FSEb}
R.~Fazio.
\newblock The {Falkner}-{Skan} equation: numerical solutions within group
  invariance theory.
\newblock {\em Calcolo}, 31:115--124, 1994.

\bibitem{Fazio:1996:NAN}
R.~Fazio.
\newblock A novel approach to the numerical solution of boundary value problems
  on infinite intervals.
\newblock {\em SIAM J. Numer. Anal.}, 33:1473--1483, 1996.

\bibitem{Fazio:1998:SAN}
R.~Fazio.
\newblock A similarity approach to the numerical solution of free boundary
  problems.
\newblock {\em SIAM Rev.}, 40:616--635, 1998.

\bibitem{Fazio:2013:BPF}
R.~Fazio.
\newblock {B}lasius problem and {F}alkner-{S}kan model: {T}{\"o}pfer's
  algorithm and its extension.
\newblock {\em Comput. \& Fluids}, 73:202--209, 2013.
\newblock A preprint related to a first version of this paper is available at
  the URL: \url{http://arxiv.org/pdf/1212.5057v1.pdf}.

\bibitem{Fazio:1990:SNA}
R.~Fazio and D.~J. Evans.
\newblock Similarity and numerical analysis for free boundary value problems.
\newblock {\em Int. J. Computer Math.}, 31:215--220, 1990.
\newblock 39~:~249, 1991.

\bibitem{Fazio:2010:ESI}
R.~Fazio and S.~Iacono.
\newblock On the equivalence of non-iterative transformation methods based on
  scaling and spiral groups.
\newblock {\em Math. Meth. Appl. Sci.}, 33:585--591, 2010.

\bibitem{Fazio:2010:TGN}
R.~Fazio and S.~Iacono.
\newblock On the translation groups and non-iterative transformation methods.
\newblock In E.~De Bernardis, R.~Spligher, and V.~Valenti, editors, {\em
  Applied and Industrial Mathematics in Italy III}, pages 331--340, Singapore,
  2010. World Scientific.

\bibitem{Fox:1969:LBL}
V.~G. Fox, L.~E. Erickson, and L.~I. Fan.
\newblock The laminar boundary layer on a moving continuous flat sheet in a
  non-newtonian fluid.
\newblock {\em AIChE J.}, 15:327--333, 1969.

\bibitem{Gallant:2012:DPS}
J.~Gallant.
\newblock {\em Doing Physics with Scientific Notebook: A Problem-solving
  Approach}.
\newblock Wiley, 2012.

\bibitem{Heun:1900:BNI}
K.~Heun.
\newblock New {Methode} zur {I}ntegration der {D}ifferentialgleichungen einer
  unbh\"{a}ngigen {V}er\"{a}nderlichen.
\newblock {\em Zeitschrift f{\"{u}}r Mathematik und Physik}, 45:23--38, 1900.

\bibitem{Kutta:1901:BNI}
M.~W. Kutta.
\newblock Beitrag zur n\"{a}herungweisen {I}ntegration totaler
  {D}ifferentialgleichungen.
\newblock {\em Zeitschrift f{\"{u}}r Mathematik und Physik}, 46:435--453, 1901.

\bibitem{Landau:HCM:1950}
H.~G. Landau.
\newblock Heat conduction in melting solid.
\newblock {\em Q. Appl. Math.}, 8:81--94, 1950.

\bibitem{Meyer:1973:IVM}
G.~H. Meyer.
\newblock {\em Initial Value Methods for Boundary Value Problems; Theory and
  Application of Invariant Imbedding}.
\newblock Academic Press, New York, 1973.

\bibitem{Na:1979:CME}
T.~Y. Na.
\newblock {\em Computational Methods in Engineering Boundary Value Problems}.
\newblock Academic Press, New York, 1979.

\bibitem{Runge:1895:UNA}
C.~Runge.
\newblock Ueber die numerische {A}ufl\"{o}sung von {D}ifferentialgleichungen.
\newblock {\em Mathematische Annalen}, 46:167--178, 1895.

\bibitem{Sachdev:1991:NOD}
P.~L. Sachdev.
\newblock {\em Nonlinear Ordinary Differential Equations and their
  Applications}.
\newblock Marcel Dekker, New York, 1991.

\bibitem{Seshadri:1985:GIE}
R.~Seshadri and T.~Y. Na.
\newblock {\em Group Invariance in Engineering Boundary Value Problems}.
\newblock Springer, New York, 1985.

\bibitem{Shampine:1997:MOS}
L.~F. Shampine and M.~W. Reichelt.
\newblock The {MATLAB ODE} suite.
\newblock {\em SIAM J. Sci. Comput.}, 18:1--22, 1997.

\bibitem{Shampine:SOM:2003}
L.F. Shampine, I.~Gladwell, and S.~Thompson.
\newblock {\em Solving {ODE}s with {MATLAB}}.
\newblock Cambridge University Press, Cambridge, 2003.

\end{thebibliography}
